\documentclass[a4paper,fleqn]{cas-sc}
\IfFileExists{stix.sty}{}{\usepackage{lmodern}\usepackage{mathptmx}}

\usepackage{amsmath, amsthm, amsfonts, amssymb, color}
\usepackage{mathrsfs}

\allowdisplaybreaks
\usepackage{pgf,tikz}
\usepackage{comment}
\usepackage[numbers]{natbib}

\newtheorem{thm}{Theorem}[section]
\newtheorem{lem}{Lemma}[section]

\newtheorem{prop}{Proposition}[section]
\theoremstyle{definition}
\newtheorem{defn}{Definition}[section]

\theoremstyle{remark}

\numberwithin{equation}{section}

\def\tsc#1{\csdef{#1}{\textsc{\lowercase{#1}}\xspace}}
\tsc{WGM}
\tsc{QE}

\begin{document}
\let\WriteBookmarks\relax
\def\floatpagepagefraction{1}
\def\textpagefraction{.001}

\shorttitle{Uncertainty principle on $\mathbb{T}^m \times \mathbb{R}^n$}

\shortauthors{Gan and Yang}

\title [mode = title]{Uncertainty principle on $\mathbb{T}^m \times \mathbb{R}^n$}

\author[1]{Wenkai Gan}
\author[1]{Yufeng Yang}

\affiliation[1]{organization={School of Mathematics, Sun Yat-sen University},
            addressline={No. 135, Xingang Xi Road, Haizhu District},
            city={Guangzhou},
            postcode={510275},
            state={Guangdong Province},
            country={China}}

\begin{abstract}
This paper establishes an uncertainty principle on the mixed space $\mathbb{T}^m \times \mathbb{R}^n$
for $n\ge 2$. Leveraging the paraboloid restriction theorem on $\mathbb{Z}^m \times \mathbb{R}^n$ derived from Barron's global-in-time Strichartz estimates for semi-periodic Schr\"odinger equation, we construct measurable sets $S \subset \mathbb{T}^m \times \mathbb{R}^n$ and $A(R,\delta) \subset \mathbb{Z}^m \times \mathbb{R}^n$ that form a strong $L^2$ annihilating pair, satisfying the uniform inequality
$\|f\|_{L^2(\mathbb{T}^m \times \mathbb{R}^n)} \leq 3 \left( \|f\|_{L^2(\mathbb{T}^m \times \mathbb{R}^n \setminus S)} + \|\hat{f}\|_{L^2(\mathbb{Z}^m \times \mathbb{R}^n \setminus A(R,\delta))} \right)$
for all $f \in L^2(\mathbb{T}^m \times \mathbb{R}^n)$, where the measures of $S$ and $A(R,\delta)$ can be made arbitrarily large while the constant $3$ remains unchanged. As an application, we derive a two-time observability inequality for the semi-periodic Schr\"odinger equation under a support condition on the phase-modified solution at the first observation time.
\end{abstract}


\maketitle

\vskip 0.5cm

\section{Introduction}
We begin by recalling the definitions of a weak (strong) $L^2$ annihilating pair. Following
  \cite{Havin}, we say that   $S, \Sigma \subset \mathbb{R}^n$ form a weak $L^2$ annihilating pair if for every $f \in L^2(\mathbb{R}^n)$ with $\operatorname{supp} f \subset S$ and $\operatorname{supp} \hat{f} \subset \Sigma$, we have $f=0$;  and  $S, \Sigma\subset \mathbb{R}^n$ form a strong $L^2$ annihilating pair if there exists $C>0$ such that for every $f \in L^2(\mathbb{R}^n)$,

$${\|f\|}_{L^2(\mathbb{R}^n)}\le C\bigl({\|f\|}_{L^2(\mathbb{R}^n\backslash S)}+{\|\hat{f}\|}_{L^2(\mathbb{R}^n \backslash \Sigma)}\bigr). $$

A fundamental problem is how to construct a pair of sets such that it is an annihilation pair.
Here we will focus on the case of $S$, $\Sigma$ being of finite Lebesgue measure. Regarding weak $L^2$ annihilation pairs, we have a fundamental fact proved by Benedicks \cite{BENEDICKS1985180} and this result was also proved by Amrein and Berthier \cite{1977On} using an alternative method:

\begin{prop}
Let $f \in L^2(\mathbb{R}^n)$. If $\operatorname{supp} f$ and $\operatorname{supp} \hat{f}$ both have finite measure in $\mathbb{R}^n$, then $f=0$.
\end{prop}

Demange and Bonami \cite{Demange} showed that it is actually a strong annihilating pair.
Therefore, for strong $L^2$ annihilation pairs, how to control the size of the constant governing the inequality is more important. Nazarov \cite{Nazarov} obtained the following inequality on $\mathbb{R}$ and Jaming \cite{JAMING200730} extended it to higher dimension:

\begin{prop}
Let $S, \Sigma$ be finite-measure sets in $\mathbb{R}^n$. There exist constants $C, C'>0$, independent of $S,\Sigma$, such that for every $f \in L^2(\mathbb{R}^n)$,
$$ \int_{\mathbb{R}^n}{|f(x)|}^2dx\le C e^{C'\min\bigl(|S||\Sigma|,\;|S|^{\frac{1}{n}}w(\Sigma),\; w(S)|\Sigma|^{\frac{1}{n}}\bigr)} \Bigl(\int_{\mathbb{R}^n \backslash S}{|f(x)|}^2dx+\int_{\mathbb{R}^n \backslash \Sigma}{|\hat{f}(\xi)|}^2d\xi\Bigr).$$
\end{prop}
\noindent
Here $w(\Sigma)$ denotes the mean width of $\Sigma$. Thus when $S$, $\Sigma$ are large, the constant in the inequality also becomes large. By contrast, Jaming, Iosevich, and Mayeli \cite{jaming2025uncertaintyprincipleannihilatingpairs} construct a pair of sets with arbitrarily large measure for which the inequality holds with a uniform constant, by applying restriction theorems on the sphere and by constructing a fattened neighborhood of the sphere:

\begin{prop}
For $n\ge 2$, there exist positive-measure sets $S', A'(R, \delta) \subset \mathbb{R}^n$ such that for every $f \in L^2(\mathbb{R}^n)$,
$$ {\|f\|}_{L^2(\mathbb{R}^n)}\le3\bigl({\|f\|}_{L^2(\mathbb{R}^n\backslash S')}+{\|\hat{f}\|}_{L^2(\mathbb{R}^n\backslash A'(R, \delta))}\bigr). $$
and the measures of $S'$ and $A'(R,\delta)$ can be made arbitrarily large while keeping the constant unchanged.
\end{prop}
\noindent

In the present work, we extend this line of investigation to the mixed space $\mathbb{T}^m \times \mathbb{R}^n$, where the
$\mathbb{T}^m$ is identified with the periodic cube ${[0, 2\pi]}^m \subset \mathbb{R}^m$. Using the tools of abstract harmonic analysis \cite{Folland} , we generalize the Fourier transform on $\mathbb{R}^n$ to $\mathbb{T}^m \times \mathbb{R}^n$.
For a function $f$ on $\mathbb{T}^m \times \mathbb{R}^n$, define its $L^p(\mathbb{T}^m \times \mathbb{R}^n)$ norm ($1\le p <\infty$) by
$${\|f\|}_{L^p(\mathbb{T}^m \times \mathbb{R}^n)}=\Bigl(\int_{\mathbb{T}^m \times \mathbb{R}^n} {|f(x, y)|}^p dxdy\Bigr)^{\frac{1}{p}}, $$
while for $p=\infty$, the norm is defined as the essential supremum of $f$.

Correspondingly, for a function $g$ on $\mathbb{Z}^m \times \mathbb{R}^n$, define its $L^p(\mathbb{Z}^m \times \mathbb{R}^n)$ norm ($1\le p < \infty$) by
$${\|g\|}_{L^p(\mathbb{Z}^m \times \mathbb{R}^n)}=\Bigl(\sum_{k \in \mathbb{Z}^m}\int_{\mathbb{R}^n}{|g(k, \xi)|}^p d\xi\Bigr)^{\frac{1}{p}}, $$
with the $L^\infty$ norm again defined as the essential supremum.

For a function $f$ on $\mathbb{T}^m \times \mathbb{R}^n$, define its Fourier transform as a function on $\mathbb{Z}^m \times \mathbb{R}^n$:
$$ \hat{f}(k, \xi)=\frac{1}{{(2\pi)}^{(m+n)/2}} \int_{\mathbb{T}^m \times \mathbb{R}^n} f(x, y) e^{-i(kx+\xi y)}\,dx\,dy. $$
For a sufficiently good $f$, the Fourier inversion formula holds:
$$f(x, y)=\sum_{k \in \mathbb{Z}^m}\frac{1}{{(2\pi)}^{(m+n)/2}}\int_{\mathbb{R}^n} \hat{f}(k, \xi)e^{i(kx+\xi y)}\,d\xi. $$

Barron \cite{Barron} obtained Strichartz estimates for the Schr\"odinger equation on $\mathbb{T}^m \times \mathbb{R}^n$, which yields a restriction theorem for the paraboloid on $\mathbb{Z}^m \times \mathbb{R}^n$. Building on this foundation, we prove our main uncertainty principle result for the mixed space:

\begin{thm}\label{thm:uncertainty}
For $n\ge 2$, there exist positive-measure sets $S \subset \mathbb{T}^m \times \mathbb{R}^n$ and $A(R, \delta) \subset \mathbb{Z}^m \times \mathbb{R}^n$ such that for every $L^2$ function $f$,
$${\|f\|}_{L^2(\mathbb{T}^m \times \mathbb{R}^n)}\le 3\Bigl({\|f\|}_{L^2(\mathbb{T}^m \times \mathbb{R}^n \backslash S)}+{\|\hat{f}\|}_{L^2(\mathbb{Z}^m\times \mathbb{R}^n \backslash A(R, \delta))}\Bigr), $$
and the measures of $S$ and $A(R,\delta)$ can be made arbitrarily large while keeping the constant unchanged.
\end{thm}

Consider the Schr\"odinger equation on $\mathbb{T}^m \times \mathbb{R}^n$:
\begin{equation}\label{eq:semiperiodic-schrodinger}
\left\{\begin{aligned}
 & i{\partial}_t u(x, y, t)+{\Delta}_{\mathbb{T}^m \times \mathbb{R}^n} u(x, y, t)=0, (x, y, t)\in \mathbb{T}^m\times \mathbb{R}^n \times (0, \infty),  \\
& u(x, y, 0)=u_0(x, y) \in L^2(\mathbb{T}^m \times \mathbb{R}^n).
\end{aligned}\right.
\end{equation}
As a corollary of Theorem~\ref{thm:uncertainty}, we derive the following two-time observability inequality:

\begin{thm}\label{thm:observability}
For $n\ge 2$, let $S$ and $A(R,\delta)$ be as in Theorem~\ref{thm:uncertainty}, and let $R'$ be the projection of $A(R,\delta)$ onto $\mathbb{R}^n$. Fix $0\le t_1<t_2$ and put $\tau=t_2-t_1$. If the solution $u$ of \eqref{eq:semiperiodic-schrodinger} satisfies
\begin{equation}\label{eq:revised-support}
\operatorname{supp}\widehat{e^{\frac{i|y|^2}{4\tau}}u(\cdot,\cdot,t_1)}
\subset A(R,\delta)\cup\bigl(\mathbb{Z}^m\times(\mathbb{R}^n\setminus R')\bigr),
\end{equation}
where the support inclusion is understood up to sets of measure zero, then
\begin{equation}\label{eq:two-time-observability}
\begin{aligned}
\int_{\mathbb{T}^m\times\mathbb{R}^n}|u_0(x,y)|^2\,dx\,dy
\le 18\Bigg(&\int_{(\mathbb{T}^m\times\mathbb{R}^n)\setminus S}|u(x,y,t_1)|^2\,dx\,dy\\
&+\int_{\mathbb{T}^m\times2\tau(\mathbb{R}^n\setminus R')}|u(x,y,t_2)|^2\,dx\,dy\Bigg).
\end{aligned}
\end{equation}
In particular, the constant is independent of $R,\delta,t_1,t_2$ and $u_0$.
\end{thm}

The remainder of this paper is structured as follows: Section 2 introduces necessary preliminary technical results including H\"older's inequality on the mixed dual space, associated embedding theorems, Barron's Strichartz estimates, and a generalized extension of Benedicks' finite-support vanishing result to the mixed space $\mathbb{T}^m \times \mathbb{R}^n$.
Section 3 is dedicated to the full proof of Theorem 1.1, via the paraboloid restriction theorem and annihilating pair construction. Section 4 presents the derivation of Theorem 1.2 as a direct application of the main uncertainty principle to semi-periodic Schr\"odinger dynamics.

\section{Preliminaries}
\label{cha:Pre thm}

In this section, we collect some foundational technical lemmas, inequalities, and generalized results required for the main theorem proofs.

We first extend the classical H\"older's inequality to the mixed discrete-continuous space $\mathbb{Z}^m \times \mathbb{R}^n$:
\begin{prop}
For any $f \in L^p(\mathbb{Z}^m \times \mathbb{R}^n), g \in L^{p'}(\mathbb{Z}^m \times \mathbb{R}^n)$, where $1\le p \le \infty, \frac{1}{p}+\frac{1}{p'}=1$, we have
$$|\sum_{k \in \mathbb{Z}^m}\int_{\mathbb{R}^n} f(k, \xi)g(k, \xi) d\xi| \le {\|f\|}_{L^p(\mathbb{Z}^m \times \mathbb{R}^n)}{ \|g\|}_{L^{p'}(\mathbb{Z}^m \times \mathbb{R}^n)}.$$
\end{prop}
\begin{proof}
First apply the H\"older inequality in $L^p(\mathbb{R}^n)$ to the integral in the $\xi$ variable to get:
$$|\sum_{k \in \mathbb{Z}^m}\int_{\mathbb{R}^n} f(k, \xi)g(k, \xi) d\xi| \le \sum_{k\in \mathbb{Z}^m} {(\int_{\mathbb{R}^n}{|f(k, \xi)|}^p d\xi)}^{\frac{1}{p}}{(\int_{\mathbb{R}^n}{|g(k, \xi)|}^{p'} d\xi)}^{\frac{1}{p'}}, $$
then apply the H\"older inequality in $l^p$ space to the $k$ variable to get:
$$\sum_{k\in \mathbb{Z}^m} {(\int_{\mathbb{R}^n}{|f(k, \xi)|}^p d\xi)}^{\frac{1}{p}}{(\int_{\mathbb{R}^n}{|g(k, \xi)|}^{p'} d\xi)}^{\frac{1}{p'}} \le {(\sum_{k \in \mathbb{Z}^m}\int_{\mathbb{R}^n}{|f(k, \xi)|}^p d\xi)}^{\frac{1}{p}}{(\sum_{k \in \mathbb{Z}^m}\int_{\mathbb{R}^n}{|g(k, \xi)|}^{p'} d\xi)}^{\frac{1}{p'}},  $$
this completes the proof.
\end{proof}

As a direct corollary of Proposition 2.1, we obtain a generalized norm embedding relation for all bounded measurable subsets of $\mathbb{Z}^m \times \mathbb{R}^n$:
\begin{prop}
Let $\Omega$ be a bounded set in $\mathbb{Z}^m \times \mathbb{R}^n$. For $1\le r < p \le \infty, f \in L^p(\Omega)$, we have
$${\|f\|}_{L^r(\Omega)} \le {\|f\|}_{L^p(\Omega)} {|\Omega|}^{\frac{1}{r}-\frac{1}{p}}. $$
\end{prop}
We now recall the multiplier framework and global-in-time Strichartz estimates originally derived by Barron, which are the core technical foundation for our paraboloid restriction theorem on the mixed domain.

Let $m(k, \xi)$ be a multiplier on $\mathbb{Z}^m \times \mathbb{R}^{n-1}$. The integral operator $T$ on $L^2(\mathbb{T}^m \times \mathbb{R}^{n-1})$ induced by $m$ is given by
$$Tf(x, y)=\sum_{k \in \mathbb{Z}^m}\frac{1}{{(2\pi)}^{(m+n-1)/2}}\int_{\mathbb{R}^{n-1}} m(k, \xi)\hat{f}(k, \xi)e^{i(kx+\xi y)}d\xi.$$
Consider the Schr\"odinger equation on $\mathbb{T}^m \times \mathbb{R}^{n-1}$:
\begin{equation}
\left\{\begin{aligned}
 & i{\partial}_t u(x, y, t)+{\Delta}_{\mathbb{T}^m \times \mathbb{R}^{n-1}} u(x, y, t)=0, (x, y, t)\in \mathbb{T}^m\times \mathbb{R}^{n-1} \times (0, \infty), \\
 & u(x, y, 0)=f(x, y) \in L^2(\mathbb{T}^m \times \mathbb{R}^{n-1}).
\end{aligned}\right.
\end{equation}
Thus $f \mapsto u$ gives an operator $Tf=e^{it{\Delta}_{\mathbb{T}^m \times \mathbb{R}^{n-1}} f(x, y)}, $ i.e.:
$$Tf(x, y, t)=\sum_{k \in \mathbb{Z}^m}\frac{1}{{(2\pi)}^{(m+n-1)/2}}\int_{\mathbb{R}^{n-1}} e^{-it({|k|}^2+{|\xi|}^2)}\hat{f}(k, \xi)e^{i(kx+\xi y)}d\xi.$$
The following Strichartz estimate is obtained by Barron \cite{Barron}:
\begin{lem}{\cite{Barron}}
Let$f\in L^2(\mathbb{T}^m \times \mathbb{R}^{n-1})$, $q \ge 2+\frac{4}{m+n-1}$, and $\alpha > (m+n-1)(\frac{1}{2}-\frac{1}{q})-\frac{2}{q}$. Then
$${\|Tf\|}_{L^q(\mathbb{T}^m \times \mathbb{R}^{n-1} \times [0, 1])}\le C {\|f\|}_{W^{\alpha, 2}(\mathbb{T}^m \times \mathbb{R}^{n-1})}, $$
where $${\|f\|}_{W^{\alpha, 2}(\mathbb{T}^m \times \mathbb{R}^{n-1})}={\|{(1+{|k|}^2+{|\xi|}^2)}^{\alpha/2}\hat{f}(k, \xi)\|}_{L^2(\mathbb{Z}^m\times \mathbb{R}^{n-1})}.$$

\end{lem}
Following the method of Benedicks \cite{BENEDICKS1985180}, we can prove Benedicks' result in the special case on $\mathbb{T}^m \times \mathbb{R}^n$ using a more specific method.
\begin{prop}
Let $f\in L^2(\mathbb{T}^m \times \mathbb{R}^n)$. If  $\operatorname{supp} f$ and $\operatorname{supp} \hat{f}$ have finite measure in $\mathbb{T}^m\times \mathbb{R}^n$ and $\mathbb{Z}^m\times \mathbb{R}^n$ respectively, then $f=0$.
\end{prop}
\begin{proof}[Proof of Proposition 2.3]

  Let $A = \operatorname{supp} f$, $B = \operatorname{supp} \hat{f}$. Without loss of generality, assume $|A| < (2\pi)^{m+n}$, otherwise we scale the $y$ variable of $f$ and consider the scaled function. Construct a function on $\mathbb{Z}^m \times \mathbb{R}^n$:
$$ \phi(k, y)=\sum_{(i, j)\in \mathbb{Z}^{m} \times \mathbb{Z}^n} \chi_{B}(k+i, y+j). $$
Let $D=(0, 0, \dots, 0) \times [0, 1]^n \subset \mathbb{Z}^m \times \mathbb{R}^n$, and consider the integral of $\phi$ over $D$:
$$
 \int_{D} \phi(k, y)=\int_{[0, 1]^n}\phi(0, y)dy=\sum_{(i, j)\in \mathbb{Z}^{m} \times \mathbb{Z}^n}\int_{[0, 1]^n} \chi_B(i, y+j)dy.
$$
Applying a translation in the $y$ variable to the integral on the far right of (4) yields
$$ \sum_{i\in \mathbb{Z}^m} \int_{\mathbb{R}^n}\chi_B(i, y)dy=|B|<\infty, $$
this shows that for almost every $y \in [0, 1]^n$, $\phi(0, y)$ is finite; i.e., for such $y$,
$ \sum_{(i, j)\in \mathbb{Z}^{m} \times \mathbb{Z}^n} \chi_B(i, y+j)$
has only finitely many non-zero terms. Consequently, for almost every $\xi \in \mathbb{R}^n$, $\hat{f}(i, \xi+j)$ is non-zero for only finitely many $(i, j)\in \mathbb{Z}^m \times \mathbb{Z}^n$.
\newline
Now fix $\xi_0 \in \mathbb{R}^n$ and construct
$$ F_{\xi_0}(x, y)=\sum_{j\in \mathbb{Z}^n}e^{-i\xi_0(y+2\pi j)}f(x, y+2\pi j), $$
then $F_{\xi_0} \in L^1(\mathbb{T}^{m+n})$, and its Fourier coefficients are:
$$ \widehat{F_{\xi_0}}(k, l)=\frac{1}{(2\pi)^{m+n}}\int_{\mathbb{T}^{m+n}}\sum_{j\in \mathbb{Z}^n}e^{-i\xi_0(y+2\pi j)}f(x, y+2\pi j)e^{-i ly-ikx}dxdy, $$
notice that by using translation in the $y$ variable and periodicity, we obtain:
$$\widehat{F_{\xi_0}}(k, l)= C \hat{f}(k, \xi_0+l).  $$
Also note that
$$\int_{\mathbb{T}^{m+n}} \sum_{j\in \mathbb{Z}^n}\chi_{A}(x, y+2\pi j) dx dy=\int_{\mathbb{T}^m \times \mathbb{R}^n} \chi_A(x, y)dxdy=|A|<(2\pi)^{n+m}, $$
this implies that the measure of the support of $F_{\xi_0}$ is less than $(2\pi)^{n+m}$, so within one period it has infinitely many zeros. However, for almost every $\xi_0 \in \mathbb{R}^n$, the Fourier coefficients of $F_{\xi_0}$ have only finitely many non-zero terms; the Fourier expansion of such $F_{\xi_0}$ is a trigonometric polynomial, hence
$$ F_{\xi_0}(x, y)=0, \text{ a.e. } \ \xi_0 \in \mathbb{R}^n, $$
and therefore
$$ \hat{f}(k, \xi_0+l)=0 , \quad \forall (k, l)\in \mathbb{Z}^m\times \mathbb{Z}^n , \text{ a.e. } \ \xi_0 \in \mathbb{R}^n. $$
By the uniqueness of the Fourier transform, $f=0$.
\end{proof}

\section{Proof of the Theorem 1.1}
In this section, via Lemma 2.1 we obtain a restriction theorem for a paraboloid on $\mathbb{Z}^m \times \mathbb{R}^n$, then thicken the paraboloid to obtain Theorem 1.1. In the next section, we derive Theorem 1.2 from Theorem 1.1.

Note that for the operator defined in Lemma 2.1, since
$$e^{i(kx+\xi y -t({|k|}^2+{|\xi|}^2))}=e^{i(x, y, t)(k, \xi, -({|k|}^2+{|\xi|}^2))}, $$
$Tf$ is essentially the Fourier transform of the function $f$ on the paraboloid
$S=\{(k, \xi, -({|k|}^2+{|\xi|}^2)) \mid k\in \mathbb{Z}^m , \xi \in \mathbb{R}^{n-1}\}$. We have the following restriction theorem:
\begin{lem}
For the paraboloid $ S=\{(k, \xi, -({|k|}^2+{|\xi|}^2))\mid  |k| \le M, |\xi|\le M, k\in \mathbb{Z}^m , \xi \in \mathbb{R}^{n-1}\} \subset \mathbb{T}^m \times \mathbb{R}^n$, let $f$ be an $L^2$ function supported on $\mathbb{T}^m \times \mathbb{R}^{n-1}\times [0, 1]$. Then
$$ {\|\hat{f}\|}_{L^2(S)}\le C {\|f\|}_{L^p(\mathbb{T}^m \times \mathbb{R}^n)}  , $$
where $M$ is a given positive number and $1\le p \le \frac{2m+2n+2}{n+m+3}. $
\end{lem}
\begin{proof}
According to Lemma 2.1, let $g\in L^2(\mathbb{T}^m\times \mathbb{R}^n)$. Then
$$
\begin{aligned}
&{\|\sum_{k \in \mathbb{Z}^m}\int_{\mathbb{R}^{n-1}} e^{-it({|k|}^2+{|\xi|}^2)}\hat{g}(k, \xi)e^{i(kx+\xi y)}d\xi\|}_{L^q(\mathbb{T}^m \times \mathbb{R}^{n-1} \times [0, 1])}
\\
  &\le C {\|g\|}_{W^{\alpha, 2}(\mathbb{T}^m \times \mathbb{R}^{n-1})} \\
&=C{\|{(1+{|k|}^2+{|\xi|}^2)}^{\alpha /2}\hat{g}(k, \xi)\|}_{L^2(\mathbb{Z}^m\times \mathbb{R}^{n-1})}.
\end{aligned}
$$
Let $h$ be an $L^2$ function supported on $\{(k, \xi)\mid |k|\le M, |\xi| \le M , k\in \mathbb{Z}^m, \xi \in \mathbb{R}^{n-1}\}$. Replace $g$ in (3.1) by $\hat h$ to obtain
    $$\begin{aligned}
&{\|\sum_{k \in \mathbb{Z}^m}\int_{\mathbb{R}^{n-1}} e^{-it({|k|}^2+{|\xi|}^2)}h(k, \xi)e^{i(kx+\xi y)}d\xi\|}_{L^q(\mathbb{T}^m \times \mathbb{R}^{n-1} \times [0, 1])} \\
&\le C {\|{(1+{|k|}^2+{|\xi|}^2)}^{\alpha /2}h(k, \xi)\|}_{L^2(\mathbb{Z}^m\times \mathbb{R}^{n-1})}\\ &\le C{\|h\|}_{L^2(\mathbb{Z}^m \times \mathbb{R}^{n-1})},
\end{aligned}$$
i.e.,
$$ {\|\widehat{h\ d\sigma_S}\|}_{L^q(\mathbb{T}^m\times \mathbb{R}^{n-1}\times [0, 1])} \le C {\|h\|}_{L^2(S)}. $$
Take $f$ to be an $L^2$ function supported on $\mathbb{T}^m \times \mathbb{R}^{n-1}\times [0, 1]$, let $\frac{1}{p}+\frac{1}{q}=1$. By H\"older's inequality,
$$\begin{aligned}
 |\int_{\mathbb{T}^m \times \mathbb{R}^n} f\widehat{h\ d\sigma_S} dxdy|&\le {\|f\|}_{L^p (\mathbb{T}^m \times \mathbb{R}^n)}{\|\widehat{h\ d\sigma_S}\|}_{L^q(\mathbb{T}^m \times \mathbb{R}^{n-1}\times [0, 1])}\\
&\le C{\|f\|}_{L^p (\mathbb{T}^m \times \mathbb{R}^n)} {\|h\|}_{L^2(S)},
\end{aligned}$$
but
$$\int_{\mathbb{T}^m \times \mathbb{R}^n} f\widehat{h\ d\sigma_S} dxdy =\int_{S} \hat{f} h d\sigma_S,  $$
by duality, we obtain
$$ {\|\hat{f}\|}_{L^2(S)}\le C {\|f\|}_{L^p(\mathbb{T}^m \times \mathbb{R}^n)} .  $$
\end{proof}
To formalize the annihilating pair concept for our mixed periodic-Euclidean space, we present the following formal definition, adapted to the dual mixed domain structure:
\begin{defn}
We call $(S, \Sigma)$ a weak annihilating pair for $\mathbb{T}^m \times \mathbb{R}^n$ and $\mathbb{Z}^m \times \mathbb{R}^n$ if for every $a \in L^2(\mathbb{T}^m \times \mathbb{R}^n)$, whenever $supp\  a \subset S$ and $supp\  \hat{a} \subset \Sigma$, we have $a=0$. We call $(S, \Sigma)$ an $L^2$ strong annihilating pair for $\mathbb{T}^m \times \mathbb{R}^n$ and $\mathbb{Z}^m \times \mathbb{R}^n$ if for every $a \in L^2(\mathbb{T}^m \times \mathbb{R}^n)$, there exists a constant $A$ such that
$$ {\|a\|}_{L^2(\mathbb{T}^m \times \mathbb{R}^n)} \le A({\|a\|}_{\mathbb{T}^m \times \mathbb{R}^n \backslash S}+{\|\hat{a}\|}_{\mathbb{Z}^m\times \mathbb{R}^n \backslash \Sigma}). $$
\end{defn}
We next derive a sufficient condition for two finite-measure sets to form such a strong annihilating pair, which will act as our core tool for constructing the desired family of sets with uniform constants:
\begin{prop}
Let $1\le p\le 2\le q $, $ S \subset \mathbb{T}^m \times \mathbb{R}^n$, $\Sigma \subset \mathbb{Z}^m \times \mathbb{R}^n$, with $S, \Sigma$ both of finite measure. If for every $f \in L^2(\mathbb{T}^m \times \mathbb{R}^n)$,
$${(\int_{\Sigma} {{|\hat{f}(k, \xi)|}^q)}^{\frac{1}{q}}}\le C_{p, q}(\Sigma){ (\int_{\mathbb{T}^m\times \mathbb{R}^n}{|f(x, y)|}^p dxdy)}^{\frac{1}{p}}, $$
and if
$$ C_{p, q}(\Sigma){|S|}^{\frac{1}{p}-\frac{1}{2}}{|\Sigma|}^{\frac{1}{2}-\frac{1}{q}}<1, $$
then $(S, \Sigma)$ is an $L^2$ strong annihilating pair, i.e.,
$$ {\|f\|}_{L^2(\mathbb{T}^m \times \mathbb{R}^n)} \le A_{ann}(S, \Sigma)({\|f\|}_{L^2(\mathbb{T}^m \times \mathbb{R}^n \backslash S)}+{\|\hat{f}\|}_{L^2(\mathbb{Z}^m\times \mathbb{R}^n \backslash \Sigma)}), $$
where
$$A_{ann}(S, \Sigma)=1+\frac{1}{1-C_{p, q}(\Sigma){|S|}^{\frac{1}{p}-\frac{1}{2}}{|\Sigma|}^{\frac{1}{2}-\frac{1}{q}}}. $$
\end{prop}
\begin{proof}
Let $\chi_S$ be the characteristic function of $S$. By Proposition 2.2,
$$
{\|\widehat{\chi_S f}\|}_{L^2(\Sigma)} \le {|\Sigma|}^{\frac{1}{2}-\frac{1}{q}}{\|\widehat{\chi_S f}\|}_{L^q(\Sigma)},
$$
by the restriction condition, applying Proposition 2.2 again, the right-hand side of (6) becomes
$$\begin{aligned}
\le {|\Sigma|}^{\frac{1}{2}-\frac{1}{q}}C_{p, q}(\Sigma){\|f\|}_{L^p(S)}&\le
C_{p, q}(\Sigma){|S|}^{\frac{1}{p}-\frac{1}{2}}{|\Sigma|}^{\frac{1}{2}-\frac{1}{q}} {\|f\|}_{L^2(S)}.
\end{aligned}$$
Moreover,
$$\begin{aligned}
{\|\widehat{ \chi_S f}\|}_{L^2(\mathbb{Z}^m\times \mathbb{R}^n \backslash \Sigma)} &\ge  {\|\widehat{ \chi_S f}\|}_{L^2(\mathbb{Z}^m\times \mathbb{R}^n) }- {\|\widehat{ \chi_S f}\|}_{L^2(\Sigma)}\\
 &={\|\chi_S f\|}_{L^2(\mathbb{T}^m \times \mathbb{R}^n) }-{\|\widehat{ \chi_S f}\|}_{L^2(\Sigma)}\\
&\ge (1-C_{p, q}(\Sigma){|S|}^{\frac{1}{p}-\frac{1}{2}}{|\Sigma|}^{\frac{1}{2}-\frac{1}{q}}){\|\chi_S f\|}_{L^2(\mathbb{T}^m \times \mathbb{R}^n)},
\end{aligned}$$
i.e.,
$$ {\|\chi_S f\|}_{L^2(\mathbb{T}^m \times \mathbb{R}^n) } \le \frac{1}{1-C_{p, q}(\Sigma){|S|}^{\frac{1}{p}-\frac{1}{2}}{|\Sigma|}^{\frac{1}{2}-\frac{1}{q}}} {\|\widehat{ \chi_S f}\|}_{L^2(\mathbb{Z}^m\times \mathbb{R}^n \backslash \Sigma)}, $$
but
$$\begin{aligned}
{\|\widehat{ \chi_S f}\|}_{L^2(\mathbb{Z}^m\times \mathbb{R}^n \backslash \Sigma)}&\le \|\hat{f}\|_{L^2(\mathbb{Z}^m\times \mathbb{R}^n \backslash \Sigma)}+\|\widehat{\chi_{\mathbb{T}^m\times \mathbb{R}^n \backslash S}f}\|_{L^2(\mathbb{Z}^m\times \mathbb{R}^n \backslash \Sigma)}\\
 &\le \|\hat{f}\|_{L^2(\mathbb{Z}^m\times \mathbb{R}^n \backslash \Sigma)}+{\|f\|}_{L^2(\mathbb{T}^m\times \mathbb{R}^n \backslash S)} ,
\end{aligned}$$
thus,
$$\begin{aligned}
{\|f\|}_{L^2(\mathbb{T}^m \times \mathbb{R}^n)} &\le {\|f\|}_{L^2(S)}+{\|f\|}_{L^2(\mathbb{T}^m\times \mathbb{R}^n \backslash S)} \\
&\le \frac{1}{1-C_{p, q}(\Sigma){|S|}^{\frac{1}{p}-\frac{1}{2}}{|\Sigma|}^{\frac{1}{2}-\frac{1}{q}}}(\|\hat{f}\|_{L^2(\mathbb{Z}^m\times \mathbb{R}^n \backslash \Sigma)}+{\|f\|}_{L^2(\mathbb{T}^m\times \mathbb{R}^n \backslash S)})+{\|f\|}_{L^2(\mathbb{T}^m\times \mathbb{R}^n \backslash S)}\\
 &\le (1+\frac{1}{1-C_{p, q}(\Sigma){|S|}^{\frac{1}{p}-\frac{1}{2}}{|\Sigma|}^{\frac{1}{2}-\frac{1}{q}}})(\|\hat{f}\|_{L^2(\mathbb{Z}^m\times \mathbb{R}^n \backslash \Sigma)}+{\|f\|}_{L^2(\mathbb{T}^m\times \mathbb{R}^n \backslash S)}).
\end{aligned}$$
This completes the proof.
\end{proof}
We establish Theorem 1.1 by applying Lemma 3.1 and thicken the paraboloid surface into a finite-measure subset $A(R, \delta)$ that satisfies the conditions of Proposition 3.1. We systematically adjust scaling parameters to make the strong annihilating pair constant exactly equal to 3, while allowing the measures of both sets to grow without bound.
\begin{proof}[Proof of Theorem 1.1]
Define
$$C=\{(k, \xi, -({|k|}^2+{|\xi|}^2))\mid  |k| \le M, |\xi|\le M, k\in \mathbb{Z}^m , \xi \in \mathbb{R}^{n-1}\}, $$
then $C$ can be viewed as a countable union of paraboloids in $\mathbb{R}^n$. Write
$$ C_k=\{(k, \xi, -({|k|}^2+{|\xi|}^2))\mid  |\xi|\le M, \xi \in \mathbb{R}^{n-1}\}, $$
we thicken each paraboloid in $C$, i.e., thicken each $C_k$. Let $\delta$ be a small positive number. Define
$$ C_{k, \delta}=\{ (k, \xi)\in \mathbb{Z}^m \times \mathbb{R}^n \mid \exists \eta \in C_k, \text{ s.t. } |\xi-\eta|< \delta\}, $$
$$C_{\delta}={\bigcup}_{k\in \mathbb{Z}^m} C_{k, \delta}, $$
then for a function $f$ supported on $\mathbb{T}^m \times \mathbb{R}^{n-1}\times [0, 1]$, applying Lemma 3.1 we obtain
$$\begin{aligned}
{(\int_{C_{\delta}} {|\hat{f}(k, \xi)|}^2)}^\frac{1}{2}&={(\sum_{|k|\le M}\int_{C_{k, \delta}} {|\hat{f}(k, \xi)|}^2)}^\frac{1}{2}\\
&={\delta}^{\frac{1}{2}}{(\sum_{|k|\le M}\int_{C_k} {|\hat{f}(k, \xi)|}^2)}^\frac{1}{2})
={\delta}^{\frac{1}{2}}{ \|\hat{f}\|}_{L^2(C)}\\ &\le {\delta}^{\frac{1}{2}}C_{p, 2}(C) {(\int_{\mathbb{T}^m\times \mathbb{R}^n} {|f(x, y)|}^p)}^{\frac{1}{p}}.
\end{aligned}$$
Let $R>0$ and define
$$ RC=\{ (k, R\xi)\mid (k, \xi)\in C\}, $$
then
$$ {(RC)}_{\delta}=R(C_{\frac{\delta}{R}}), $$
define
$$ f_R(x, y)=\frac{1}{R^n}f(x, \frac{y}{R}) , $$
so that
$$\hat{f_R}(k, \xi)=\hat{f}(k, R\xi), $$
hence
$$\begin{aligned}
 {(\int_{{(RC)}_{\delta}} {|\hat{f}(k, \xi)|}^2)}^{\frac{1}{2}}&={(\int_{R(C_{\frac{\delta}{R}})} {|\hat{f}(k, \xi)|}^2)}^{\frac{1}{2}} \\ &=R^{\frac{n}{2}}{(\int_{C_{\frac{\delta}{R}}} {|\hat{f}(k, R\eta)|}^2)}^{\frac{1}{2}}\\
&=R^{\frac{n}{2}}{(\int_{C_{\frac{\delta}{R}}} {|\hat{f_R}(k, \eta)|}^2)}^\frac{1}{2}\le R^{\frac{n}{2}}{(\frac{\delta}{R})}^{\frac{1}{2}}C_{p, 2}(C){(\int_{\mathbb{T}^m\times \mathbb{R}^n} {|f_R(x, y)|}^p)}^{\frac{1}{p}}\\
&=C_{p, 2}(C)R^{n(\frac{1}{p}+\frac{1}{2}-1)}{(\frac{\delta}{R})}^{\frac{1}{2}}{(\int_{\mathbb{T}^m\times \mathbb{R}^n} {|f(x, y)|}^p)}^{\frac{1}{p}}.
\end{aligned}$$
Now we assume that f supported on $\mathbb{T}^m \times \mathbb{R}^{n-1}\times [0, N]$, $N>0$.
We define
$$ f_{N}(x,y)=f(x,Ny), $$
then
$$ \hat{f_N}(k, \xi)=\frac{1}{N^n}\hat{f}(k, \frac{\xi}{N}), $$
so
$$\begin{aligned}
 {(\int_{{(RC)}_{\delta}} {|\hat{f_N}(k, \xi)|}^2)}^{\frac{1}{2}}
 &=\frac{1}{N^{n}}{(\int_{R(C_{\frac{\delta}{R}})} {|\hat{f}(k, \frac{\xi}{N})|}^2)}^{\frac{1}{2}} \\
&=\frac{1}{N^{\frac{n}{2}}}{(\int_{\frac{R}{N}(C_{\frac{\delta}{R}})} {|\hat{f}(k, \xi)|}^2)}^{\frac{1}{2}}. \\
\end{aligned}$$
Using the inequality above,
$$\begin{aligned}
 \frac{1}{N^{\frac{n}{2}}}{(\int_{\frac{R}{N}(C_{\frac{\delta}{R}})} {|\hat{f}(k, \xi)|}^2)}^{\frac{1}{2}}
&\leqslant C_{p, 2}(C)R^{n(\frac{1}{p}+\frac{1}{2}-1)}{(\frac{\delta}{R})}^{\frac{1}{2}}{(\int_{\mathbb{T}^m\times \mathbb{R}^n} {|f_{N}(x, y)|}^p)}^{\frac{1}{p}}\\
&\leqslant C_{p, 2}(C)R^{n(\frac{1}{p}+\frac{1}{2}-1)}{(\frac{\delta}{R})}^{\frac{1}{2}}{N}^{-\frac{n}{p}}{(\int_{\mathbb{T}^m\times \mathbb{R}^n} {|f(x, y)|}^p)}^{\frac{1}{p}},
\end{aligned}$$
replace $R$ with $NR$,$\delta$ with $N\delta$, we have
$$\begin{aligned}
{(\int_{R(C_{\frac{\delta}{R}})} {|\hat{f}(k, \xi)|}^2)}^{\frac{1}{2}}
&\leqslant C_{p, 2}(C)R^{n(\frac{1}{p}+\frac{1}{2}-1)}{(\frac{\delta}{R})}^{\frac{1}{2}}{N}^{n(\frac{1}{p}-\frac{1}{2})}{N}^{n(\frac{n}{2}-\frac{n}{p})}{(\int_{\mathbb{T}^m\times \mathbb{R}^n} {|f(x, y)|}^p)}^{\frac{1}{p}}\\
&\leqslant C_{p, 2}(C)R^{n(\frac{1}{p}-\frac{1}{2})}{(\frac{\delta}{R})}^{\frac{1}{2}}{(\int_{\mathbb{T}^m\times \mathbb{R}^n} {|f(x, y)|}^p)}^{\frac{1}{p}},
\end{aligned}$$
use the same method,we can easily prove that applying a translation to such a function does not affect the result.

Now consider the general case.Give a function $f$,we take two sequence
$$ 0=g_{0}<g_{1}<\cdots<g_{n}<\cdots, \lim_{k\rightarrow \infty }g_{k}= \infty, $$
$$ 0=h_{0}>h_{1}>\cdots>h_{n}>\cdots, \lim_{k\rightarrow \infty }h_{k}= -\infty, $$
such that
$$ {\|f^{+}_{k}\|}_{p}\leqslant \frac{1}{2^k}{\|f\|}_{p}, $$
$$ {\|f^{-}_{k}\|}_{p}\leqslant \frac{1}{2^k}{\|f\|}_{p}, $$
where
$$ f^{+}_{k}=\chi_{A_{k}}f, $$
$$ f^{-}_{k}=\chi_{B_{k}}f, $$
$$ A_{k}=\mathbb{T}^m \times \mathbb{R}^{n-1}\times [g_{k}, g_{k+1}], $$
$$ B_{k}=\mathbb{T}^m \times \mathbb{R}^{n-1}\times [h_{k}, h_{k+1}], $$
so
$$ f=\sum_{k=0}^{\infty }f^{+}_{k}+ \sum_{k=0}^{\infty }f^{-}_{k}. $$
From the above discussion, we have
$$\begin{aligned}
{(\int_{R(C_{\frac{\delta}{R}})} {|\hat{f}(k, \xi)|}^2)}^{\frac{1}{2}}
&\leqslant \sum_{k=0}^{\infty }{\Big({(\int_{R(C_{\frac{\delta}{R}})} {|\hat{f^{+}_{k}}(k, \xi)|}^2)}^{\frac{1}{2}}
+{(\int_{R(C_{\frac{\delta}{R}})} {|\hat{f^{-}_{k}}(k, \xi)|}^2)}^{\frac{1}{2}}\Big)}\\
&\leqslant \sum_{k=0}^{\infty }{C_{p, 2}(C)R^{n(\frac{1}{p}-\frac{1}{2})}{(\frac{\delta}{R})}^{\frac{1}{2}}
\Big({(\int_{\mathbb{T}^m\times \mathbb{R}^n} {|f^{+}_{k}(x, y)|}^p)}^{\frac{1}{p}}+
{(\int_{\mathbb{T}^m\times \mathbb{R}^n} {|f^{+}_{k}(x, y)|}^p)}^{\frac{1}{p}} \Big) }\\
&\leqslant \sum_{k=0}^{\infty }{C_{p, 2}(C)R^{n(\frac{1}{p}-\frac{1}{2})}{(\frac{\delta}{R})}^{\frac{1}{2}}
\frac{1}{2^{k-1}}{(\int_{\mathbb{T}^m\times \mathbb{R}^n} {|f(x, y)|}^p)}^{\frac{1}{p}} }\\
&\leqslant 4{C_{p, 2}(C)R^{n(\frac{1}{p}-\frac{1}{2})}{(\frac{\delta}{R})}^{\frac{1}{2}}
{(\int_{\mathbb{T}^m\times \mathbb{R}^n} {|f(x, y)|}^p)}^{\frac{1}{p}} }. \\
\end{aligned}$$

Set
$$A(R, \delta)={(RC)}_{\delta}, $$
then for sufficiently large $R$ and sufficiently small $\delta$,
$$|A(R, \delta)|\approx R^{n-1} \delta , $$
take $p=\frac{4}{3}$, $C'_{p, 2}(A(R, \delta))=4C_{p, 2}(C)\delta^\frac{1}{2} R^{\frac{n}{4}-\frac{1}{2}}$, we have
$${(\int_{A(R, \delta)} {{|\hat{f}(k, \xi)|}^2)}^{\frac{1}{2}}}\le C'_{p, 2}(A(R, \delta)){ (\int_{\mathbb{T}^m\times \mathbb{R}^n}{|f(x, y)|}^{p} dxdy)}^{\frac{1}{p}}. $$

For a finite-measure set $S \subset \mathbb{T}^m \times \mathbb{R}^n$, choose $\delta, R$ such that
$$ 4C_{p_M, 2}(C)\delta^\frac{1}{2} R^{\frac{n}{4}-\frac{1}{2}}{|S|}^{\frac{1}{4}}<1, $$
by Proposition 3.1, we obtain
$$ {\|f\|}_{L^2(\mathbb{T}^m \times \mathbb{R}^n)} \le A_{ann}(S, A(R, \delta))({\|f\|}_{L^2(\mathbb{T}^m \times \mathbb{R}^n \backslash S)}+{\|\hat{f}\|}_{L^2(\mathbb{Z}^m\times \mathbb{R}^n \backslash A(R, \delta))}), $$
where
$$A_{ann}(S, A(R, \delta))=1+\frac{1}{1-4C_{p, 2}(C)\delta^\frac{1}{2} R^{\frac{n}{4}-\frac{1}{2}}{|S|}^{\frac{1}{4}}}, $$
take
$$\delta={(8C_{p_M, 2}(C) R^{\frac{n}{4}-\frac{1}{2}}{|S|}^{\frac{1}{4}})}^{-2}, $$
then
$$ |A(R, \delta)| \approx  R^{\frac{n}{2}} {|S|}^{-\frac{1}{2}},  $$
and
$$A_{ann}(S, A(R, \delta))=3. $$

Note that since $\delta$ is determined by $R$ and $|S|$, we can first take $|S|$ very large and then $R$ sufficiently large, so that $|A(R, \delta)|$ can also be arbitrarily large. That is, when the constant in the inequality is fixed at $3$, the measures of $S$ and $A(R, \delta)$ can both be arbitrarily large.
\end{proof}

\section{Proof of Theorem 1.2}
We now establish Theorem~\ref{thm:observability} as an application of Theorem~\ref{thm:uncertainty}, keeping the normalization constants in the Fourier transform. The relation between uncertainty principles and two-time observability for the free Schr\"odinger equation is also discussed by Wang, Wang and Zhang \cite{WangWangZhang2019}.

\begin{proof}[Proof of Theorem~\ref{thm:observability}]
Write $G=\mathbb{T}^m\times\mathbb{R}^n$, $\widehat G=\mathbb{Z}^m\times\mathbb{R}^n$, $A=A(R,\delta)$ and $B=\mathbb{Z}^m\times R'$. Set
$$
\tau=t_2-t_1,\qquad w(x,y)=u(x,y,t_1),\qquad v(x,y)=e^{\frac{i|y|^2}{4\tau}}w(x,y).
$$
Since the Schr\"odinger group is unitary on $L^2(G)$ and the phase factor has modulus one,
\begin{equation}\label{eq:phase-norms}
\|v\|_{L^2(G)}=\|w\|_{L^2(G)}=\|u_0\|_{L^2(G)},
\qquad
\|v\|_{L^2(G\setminus S)}=\|u(t_1)\|_{L^2(G\setminus S)}.
\end{equation}
Taking the Fourier transform of \eqref{eq:semiperiodic-schrodinger} and solving the resulting ordinary differential equation from time $t_1$ gives
$$
\widehat u(k,\xi,t_2)=e^{-i\tau(|k|^2+|\xi|^2)}\widehat w(k,\xi).
$$
For $k\in\mathbb{Z}^m$, define
$$
w_k(y)=\frac{1}{(2\pi)^{m/2}}\int_{\mathbb{T}^m}w(x,y)e^{-ik\cdot x}\,dx,
$$
and define $v_k$ and $u_k(\cdot,t_2)$ in the same way. Then $v_k(y)=e^{i|y|^2/(4\tau)}w_k(y)$ and, with the normalization $(2\pi)^{-n/2}$ for the Fourier transform on $\mathbb{R}^n$,
$$
\widehat{v_k}(\xi)=\widehat v(k,\xi),\qquad
u_k(\cdot,t_2)=e^{-i\tau|k|^2}e^{i\tau\Delta_{\mathbb{R}^n}}w_k.
$$
The free Schr\"odinger kernel therefore yields
\begin{equation}\label{eq:mode-kernel}
\begin{aligned}
u_k(y,t_2)
&=\frac{e^{-i\tau|k|^2}}{(4\pi i\tau)^{n/2}}
\int_{\mathbb{R}^n}e^{\frac{i|y-z|^2}{4\tau}}w_k(z)\,dz\\
&=\frac{e^{-i\tau|k|^2+\frac{i|y|^2}{4\tau}}}{(4\pi i\tau)^{n/2}}
\int_{\mathbb{R}^n}e^{-\frac{iy\cdot z}{2\tau}}v_k(z)\,dz\\
&=(2i\tau)^{-n/2}e^{-i\tau|k|^2+\frac{i|y|^2}{4\tau}}
\widehat v\left(k,\frac{y}{2\tau}\right).
\end{aligned}
\end{equation}
Here $(2i\tau)^{-n/2}=(2\tau)^{-n/2}e^{-i\pi n/4}$. The calculation is first made for Schwartz functions. Since the Schr\"odinger propagator, the Fourier transform, multiplication by the phase factors and the normalized dilation are unitary on $L^2(\mathbb{R}^n)$, the first and last expressions in \eqref{eq:mode-kernel} agree in the $L^2$ sense for every $w_k\in L^2(\mathbb{R}^n)$ by density. Thus the formula applies to the given $L^2$ solution without imposing the support hypothesis on any approximating sequence.

Substituting $v$ into Theorem~\ref{thm:uncertainty}, using $(a+b)^2\le 2(a^2+b^2)$ and \eqref{eq:phase-norms}, we obtain
\begin{equation}\label{eq:uncertainty-applied}
\|u_0\|_{L^2(G)}^2
\le 18\left(\|u(t_1)\|_{L^2(G\setminus S)}^2
+\|\widehat v\|_{L^2(\widehat G\setminus A)}^2\right).
\end{equation}
It remains to express the second term as an observation at time $t_2$. Since $R'$ is the projection of $A$, we have $A\subset B$ and
$$
\widehat G\setminus A=(\widehat G\setminus B)\,\dot\cup\,(B\setminus A).
$$
The support condition \eqref{eq:revised-support} is equivalent to $\widehat v=0$ almost everywhere on $B\setminus A$. Consequently,
\begin{equation}\label{eq:frequency-complement}
\begin{aligned}
\|\widehat v\|_{L^2(\widehat G\setminus A)}^2
&=\|\widehat v\|_{L^2(\widehat G\setminus B)}^2
+\|\widehat v\|_{L^2(B\setminus A)}^2\\
&=\sum_{k\in\mathbb{Z}^m}\int_{\mathbb{R}^n\setminus R'}|\widehat v(k,\xi)|^2\,d\xi.
\end{aligned}
\end{equation}
On the other hand, integrating first in $x$, using Parseval's identity on $\mathbb{T}^m$ and then applying \eqref{eq:mode-kernel}, we have
\begin{equation}\label{eq:observation-identity}
\begin{aligned}
&\int_{\mathbb{T}^m\times2\tau(\mathbb{R}^n\setminus R')}|u(x,y,t_2)|^2\,dx\,dy\\
&\qquad=\sum_{k\in\mathbb{Z}^m}\int_{2\tau(\mathbb{R}^n\setminus R')}|u_k(y,t_2)|^2\,dy\\
&\qquad=(2\tau)^{-n}\sum_{k\in\mathbb{Z}^m}\int_{2\tau(\mathbb{R}^n\setminus R')}
\left|\widehat v\left(k,\frac{y}{2\tau}\right)\right|^2\,dy\\
&\qquad=\sum_{k\in\mathbb{Z}^m}\int_{\mathbb{R}^n\setminus R'}|\widehat v(k,\xi)|^2\,d\xi.
\end{aligned}
\end{equation}
The last equality follows from the change of variables $y=2\tau\xi$; the interchange of the nonnegative sums and integrals is justified by Tonelli's theorem. Combining \eqref{eq:frequency-complement} and \eqref{eq:observation-identity}, we conclude that
$$
\|\widehat v\|_{L^2(\widehat G\setminus A)}^2
=\int_{\mathbb{T}^m\times2\tau(\mathbb{R}^n\setminus R')}|u(x,y,t_2)|^2\,dx\,dy.
$$
Substituting this identity into \eqref{eq:uncertainty-applied} proves \eqref{eq:two-time-observability}. Since $\tau=t_2-t_1$ and $w=u(t_1)$ were used from the outset, the support condition is imposed at the first observation time throughout. The proof is complete.
\end{proof}

\printcredits

\bibliographystyle{unsrtnat}
\bibliography{refs}

\end{document}